\documentclass{amsart}
\usepackage{amsmath, amsthm, amssymb, amscd, epsfig}

\input xy
\xyoption{all}

\newcommand{\Bruhatl}{l_{\emptyset}}
\newcommand{\Bruhato}{\leq_{\emptyset}}

\newcommand{\supp}{\operatorname{Supp}}
\newcommand{\twistl}[1]{l_{#1}}
\newcommand{\twisto}[1]{\leq_{#1}}
\newcommand{\op}{{\text{\rm op}}}

\swapnumbers
\theoremstyle{plain}  
\newtheorem{theorem}[subsection]{Theorem}
\newtheorem{proposition}[subsection]{Proposition}
\newtheorem{lemma}[subsection]{Lemma}
\newtheorem{corollary}[subsection]{Corollary}

\theoremstyle{definition}

\theoremstyle{remark}
\newtheorem{remark}[subsection]{Remark}

\begin{document}

\title{On the combinatorics of $B\times B$-orbits 
                        on group compactifications}

\author{Yu Chen} \author{Matthew Dyer}
\address{Department of Mathematics \\ University of Notre Dame \\
         Notre Dame, Indiana 46556, U.S.A.}
\date{}

\begin{abstract}
It is shown that there is an order isomorphism $\varphi'$ from the
poset $V$ of $B\times B$-orbits on the wonderful compactification of a
semi-simple adjoint group $G$ with Weyl group $W$ to an interval in
reverse Chevalley-Bruhat order on a non-canonically associated
Coxeter group $\hat{W}$ (in general neither finite nor affine).
Moreover, $\varphi'$ preserves the corresponding Kazhdan-Lusztig
polynomials.
Springer's (partly conjectural) construction of Kazhdan-Lusztig
polynomials for the analogues of $V$ for general Coxeter groups $W$ is
completed by reducing it by a similar order isomorphism to known
results involving a ``twisted'' Chevalley-Bruhat order on $\hat{W}$.
\end{abstract}

\maketitle


\section*{Introduction}

Let $G$ be a connected semi-simple adjoint algebraic group over an
algebraically closed field.
Denote by $B$ a Borel subgroup of $G$ and by $T$ a maximal torus in $B$.
$T$ and $B$ determine the Weyl group $W$ of $G$ and the set $S$ of its
simple reflections.
For any $I \subseteq S$, we write $W_I$ for the standard parabolic 
subgroup of $W$ generated by $I$ and $W^I$ for the set of minimal left
coset representatives of $W_I$ in $W$.

According to \cite[\S 1]{Springer:intersection},
$G$ has a ``wonderful'' compactification $X$,
which is a smooth irreducible projective $G \times G$-variety and 
contains $G$ as an open subvariety.
The $B \times B$-orbits on $X$ can be parametrized as $\mathcal{O}_v$
for $v\in V$, 
where $V=\{\,[I,a,b] \mid I \subseteq S, a \in W^I, b \in W \,\}$.
$V$ is a poset endowed with the partial order $w \leq v$ if
$\mathcal{O}_w$ is contained in the closure of $\mathcal{O}_v$ in $X$.
By \cite[Proposition 2.4]{Springer:intersection},
$[I,a_1,b_1] \leq [J,a_2,b_2]$ in $V$ if and only if
$I \subseteq J$ and there exist $c \in W_I$, $d \in W_J \cap W^I$
with $a_2 d c^{-1} \Bruhato a_1$, $b_1 c \Bruhato b_2 d$, and
$\Bruhatl (b_2 d)=\Bruhatl (b_2)+\Bruhatl (d)$,
where $\Bruhato$ and $\Bruhatl$ are the Chevalley-Bruhat order and 
the length function on $W$.

Let $\mathcal{R}=\mathbb{Z}[u,u^{-1}]$ be the integral Laurent polynomial 
ring in the indeterminate $u$,
$\mathcal{M}$ be the free left $\mathcal{R}$-module with a basis
$\{\,\tilde m_v \mid v \in V \,\}$,
and $\mathcal{H}$ be the Iwahori-Hecke algebra of $W$.
By \cite[Lemma 3.2]{Springer:intersection},
$\mathcal{M}$ admits an $(\mathcal{H},\mathcal{H})$-bimodule structure.
Springer constructs a semi-linear 
$(\mathcal{H}, \mathcal{H})$-bimodule map
$\Delta\colon \mathcal{M} \rightarrow \mathcal{M}$ in
\cite[3.3]{Springer:intersection} with respect to a certain involution on
$\mathcal{H} \otimes_{\mathcal{R}} \mathcal{H}$.
Write $\Delta (\tilde m_v)=\sum_{w \in V} \tilde b_{w,v} \tilde m_w$
for some $\tilde b_{w,v} \in \mathcal R$ ($v,w \in V$).
$\Delta$ is an involution, that is, $\Delta^2=1$
(see \cite[\S 3]{Springer:intersection}) 
and the $\tilde b_{w,v}$ may be regarded as (appropriately normalized)
analogues of the Kazhdan-Lusztig $R$-polynomials $R_{x,y}$.
Using the $\tilde b_{w,v}$ in place  of the $R_{x,y}$, Springer
defines analogues $c_{w,v}$ of the Kazhdan-Lusztig polynomials
$P_{x,y}$ and shows that they compute the Poincar\'e series
at a point in $\mathcal{O}_w$ of the local intersection
cohomology of the closure of $\mathcal{O}_v$.

Now let $(\hat{W},\hat{S})$ be any Coxeter system containing $(W,S)$ as
a standard parabolic subgroup and such that there is a bijection
$\theta\colon S\rightarrow R:=\hat{S}\smallsetminus S$ 
such that for $r,s\in S$, $r$ and $\theta(s)$ are joined by an edge 
in the Coxeter graph of $(\hat{W},\hat{S})$ 
(i.e. do not commute in $\hat{W}$) iff $r=s$. 
Let $c_R$ be a Coxeter element of the standard parabolic subgroup 
$\hat{W}_R$ and let $w_S$ be the longest element of $W=\hat{W}_S$. 
Let $\Omega$ denote the interval $[w_Sc_Rw_S,1]$ in $\hat{W}$ in the
order induced by the reverse Chevalley-Bruhat order. 
The main result of this paper implies that there is a poset isomorphism
$\varphi'\colon V \rightarrow \Omega$ satisfying
$c_{w,v}=Q_{\varphi'(v),\varphi'(w)}$,  where the $Q_{x,y}$ denote the
inverse Kazhdan-Lusztig polynomials \cite[(2.1.6)]{Kazhdan:Schubert}
for $\hat{W}$.

More generally, let $(W,S)$ be a Coxeter system with finitely
many simple reflections. 
Springer extends in \cite[\S 6]{Springer:intersection} the 
constructions of the poset $V$, $(\mathcal{H},\mathcal{H})$-bimodule
$\mathcal M$, and map $\Delta$ to $(W,S)$; 
the only needed modification in this more general setting is thinking
of $\Delta$ as a semi-linear map from $\mathcal{M}$ to its completion
$\tilde{\mathcal{M}}$.
Springer conjectures in \cite[6.10]{Springer:intersection} that
$\Delta$ is still an involution in the general setting, which
would imply that there exist analogues of the Kazhdan-Lusztig
polynomials $c_{w,v}$ and inverse Kazhdan-Lusztig polynomials 
$c_{w,v}^{\mbox{\tiny inv}}$ for $V$ 
(see \cite[Formulas (7) and (8)]{Kallen:computing}).

On the other hand, in \cite{Dyer:Iwahori-Hecke,Dyer:Hecke2}, there are
attached to any suitable subset $A$ (initial section of a
reflection order) of the reflections $\hat{T}$ of an arbitrary
Coxeter system $(\hat{W},\hat{S})$ the analogue $\twisto{A}$ of
Chevalley-Bruhat order, analogues $R^A_{x,y}$ of the $R$-polynomials
$R_{x,y}$, analogues $P_A(x,y)$ of the Kazhdan-Lusztig polynomials
$P_{x,y}$ etc;  
to define $\twisto{A}$, $R^A_{x,y}$, $P_A(x,y)$ etc, one just regards the
standard length function on which suitable standard definitions 
(of Chevalley-Bruhat order, $R_{x,y}$, $P_{x,y}$ etc)
implicitly depend as a parameter and replaces it by a more general
length function $\twistl{A}$ depending on $A$ (the fact one obtains
well-defined notions for general $\twistl{A}$ is not obvious).

In this paper, we prove
Springer's conjecture by showing that for $(\hat{W},\hat{S})$
associated to $(W,S)$ exactly as for finite Weyl groups and 
$A:=\hat{T} \smallsetminus W$, Springer's $V$, $\mathcal{M}$,
$\Delta$, $\tilde b_{w,v}$, $c_{w,v}$ etc may be described
directly in terms of corresponding objects attached to $\twisto{A}$
on $\hat{W}$. 
In the case of a finite Coxeter group $W$, $\hat{W}$ in $\twisto{A}$ 
is order isomorphic to $\hat{W}$ in reverse Chevalley-Bruhat order, 
and we recover the above-mentioned interpretation of $c_{w,v}$ as an
inverse Kazhdan-Lusztig polynomial.

The arrangement of this paper is as follows.
Section \ref{section:Springer's conjecture} briefly recalls
some details of Springer's combinatorial constructions.
Section \ref{section:results} describes our results and lists some of
their immediate consequences. 
Section \ref{section:twisted Chevalley-Bruhat order} collects the
general properties of orders $\twisto{A}$ we shall require; 
only a few of these are not already explicit in
\cite{Dyer:Iwahori-Hecke}. 
Section \ref{section:special order} gives results we need
concerning an arbitrary Coxeter system $(\hat{W},\hat{S})$ in the
specific order $\twisto{A}$,
where $A=\hat{T} \smallsetminus W$ and $(W,S)$ is a standard  parabolic
subsystem of $(\hat{W},\hat{S})$. 
The main result there is an identity 
(Lemma \ref{lemma:recursive formula for R^A}) expressing
$R^A_{azb,1}$ for $a,b\in W$ and 
$z \in \hat{W}_{\hat{S}\smallsetminus S}$ in terms of the classical
$R$-polynomials of $(\hat{W},\hat{S})$. 
In Section \ref{section:proofs}, we finally  prove the
results of Section \ref{section:results} by showing that for
$(\hat{W},\hat{S})$ associated to $(W,S)$ as described previously,
this identity specializes to the initial condition
\ref{subsection:the map Delta} (b) for Springer's recurrence formula
for $\tilde b_{w,v}$.

We thank Professor Springer for his helpful comments on an early 
version of this manuscript.

\section{Springer's Combinatorial Constructions}
\label{section:Springer's conjecture}

Let $(W,S)$ be a Coxeter system with finitely many simple reflections.
Denote the length function on $W$ by $\Bruhatl$, the Chevalley-Bruhat 
order on $W$ by $\Bruhato$, and the identity of $W$ by $1$.
If $I \subseteq S$, then write $W_I$ for the standard parabolic 
subgroup of $W$ generated by $I$.
$W^I=\{\,x \in W \mid x \Bruhato xs \mbox{ for all }s \in I \,\}$
denotes the set of minimal left coset representatives of $W_I$ in $W$.
It is well-known that for a fixed $I \subseteq S$,
$\Bruhatl(xy)=\Bruhatl(x)+\Bruhatl(y)$ for $x\in W^I$ and $y\in  W_I$,
and that every element of $W$ is uniquely expressible in the form
$xy$ for some such $x$ and $y$.

\subsection{The poset $(V,\leq)$}
In his paper \cite[6.4]{Springer:intersection},
Springer introduces a poset
\[ V=\{\,[I,a,b] \mid I \subseteq S, a \in W^I, b \in W \,\}. \]
Let $[I,a_1,b_1]$, $[J,a_2,b_2] \in V$.
Define $[I,a_1,b_1] \leq_1 [J,a_2,b_2]$
if $I \subseteq J$, $a_2 \Bruhato a_1$, $b_1 \Bruhato b_2$,
and define $[I,a_1,b_1] \leq_2 [J,a_2,b_2]$ if
$I \subseteq J$, $b_1=b_2c$ for some $c \in W_J$ with
$\Bruhatl(b_1)=\Bruhatl(b_2)+\Bruhatl(c)$, $a_2c \Bruhato a_1$.
The partial order $\leq$ on $V$ is generated by $\leq_1$ and $\leq_2$.
That is, $v \leq w$ iff there is a sequence
$v=v_0, v_1, \ldots, v_{n-1}, v_{n}=w$ in $V$ such that
either $v_{i-1} \leq_1 v_i$ or $v_{i-1} \leq_2 v_i$ for 
$i=1, 2, \ldots, n$.

Given $[I,a,b] \in V$, define
\[ d([I,a,b])=-\Bruhatl(a)+\Bruhatl(b)+|I|, \]
where $|I|$ is the cardinality of $I$.
Springer calls $d$ the dimension function on $V$ due to its relation
to the dimensions of $B \times B$-orbits on $X$ in the case $W$ is a
finite Weyl group (see \cite[Lemma 1.3]{Springer:intersection}).

\subsection{The Iwahori-Hecke algebra of $(\mathbf{W},\mathbf{S})$}
Set $\mathbf{W}=W \times W$ and
$\mathbf{S}=(S \times \{1\}) \cup (\{1\} \times S)$.
$(\mathbf{W},\mathbf{S})$ is a Coxeter system.
By abuse of notation, we write $\Bruhato$ and $\Bruhatl$ for the 
Chevalley-Bruhat order and length function on $\mathbf{W}$,
respectively.

Let $\mathcal{R}=\mathbb{Z}[u,u^{-1}]$ be the integral Laurent 
polynomial ring in the indeterminate $u$.
$\mathcal{R}$ admits a ring involution $a \mapsto \bar a$ satisfying
$\bar u=u^{-1}$.
To simplify our notation, we define $\alpha=u-u^{-1}$.
$\alpha$ satisfies the equation $\bar \alpha=-\alpha$.
Note that $u$ corresponds to the notation $q^{1/2}$ used in 
\cite[\S 1]{Kazhdan:representations}.

The Iwahori-Hecke algebra $\mathbf{H}$ of $(\mathbf{W},\mathbf{S})$ is
the free left $\mathcal{R}$-module with a basis
$\{\,\tilde T_x \mid x \in \mathbf{W}\,\}$ and subject to the
multiplication law
\begin{equation*}
\tilde T_{\mathbf{s}} \tilde T_x=\left\{
\begin{array}{ll}
\tilde T_{{\mathbf{s}}x}                
                         &\mbox{ if } x \Bruhato {\mathbf{s}}x \\
\tilde T_{{\mathbf{s}}x}+\alpha \tilde T_x
                         &\mbox{ if } {\mathbf{s}} x \Bruhato x,
\end{array} \right.
\end{equation*}
where ${\mathbf{s}} \in \mathbf{S}$ and $x \in \mathbf{W}$.
$\tilde{T}_1$ is the multiplicative identity of $\mathbf{H}$ and every
$\tilde{T}_x$ is invertible in $\mathbf{H}$.

\subsection{The $\mathbf{H}$-module $\mathcal{M}$}
\label{subsection:an H-module}
If $a \in W^I$ and $s \in S$, then there are exactly three
possibilities for $sa$:
\begin{enumerate}
\item[(a)]
$a \Bruhato sa$ and $sa \in W^I$;
\item[(b)]
$a \Bruhato sa$ and $sa=at$ for some $t \in I$;
\item[(c)]
$sa \Bruhato a$ in which case $sa \in W^I$.
\end{enumerate}
By \cite[6.4 (b)]{Springer:intersection} there is a left action of 
$\mathbf W$ on
$V$ satisfying
\begin{equation*}
(s,1).[I,a,b]=\left\{
\begin{array}{ll}
[I,sa,b]   &\mbox{ if } s, a \mbox{ satisfy (a) or (c)}  \cr
[I,a,bt]   &\mbox{ if } s, a \mbox{ satisfy (b)},
\end{array} \right.
\end{equation*}
\begin{equation*}
(1,s).[I,a,b]=[I,a,sb]
\end{equation*}
for $s \in S$ and $[I,a,b] \in V$.

Let $\mathcal{M}$ be the free left $\mathcal{R}$-module with a basis
$\{ \tilde m_v \mid v \in V \}$.
By \cite[6.3]{Springer:intersection} there is an $\mathbf{H}$-module
structure on $\mathcal{M}$ defined by
\begin{equation*}
\tilde T_{\mathbf{s}}.\tilde{m}_v=\left\{
\begin{array}{ll}
\tilde m_{{\mathbf{s}}.v}            
               &\mbox{ if } d(v)<d({\mathbf{s}}.v) \\
\tilde m_{{\mathbf{s}}.v}+\alpha \tilde m_{v}
               &\mbox{ if } d({\mathbf{s}}.v)<d(v)
\end{array} \right.
\end{equation*}
for $\mathbf{s} \in \mathbf{S}$.
Note that $\tilde m_{v}=u^{-d(v)} m_v$ for the $m_v$ appearing in
\cite[3.1]{Springer:intersection}.

\subsection{The map $\Delta$} \label{subsection:the map Delta}
Let $\tilde{\mathcal{M}}$ be the completion of $\mathcal{M}$ consisting
of formal $\mathcal{R}$-linear combinations of elements $\tilde m_v$ 
for $v \in V$. 
Springer defines a semi-linear map 
$\Delta\colon \mathcal{M} \rightarrow \tilde{\mathcal{M}}$ 
with respect to the Kazhdan-Lusztig involution on
$\mathcal{H}\otimes_\mathcal{R}\mathcal{H}$  in
\cite[6.8]{Springer:intersection}. 
Write
\[ \Delta (\tilde m_v)=\sum_{w \in V} \tilde b_{w,v} \tilde m_w \]
for some $\tilde b_{w,v} \in \mathcal{R}$ ($w,v \in V$).
By \cite[6.9(i), 3.8, 6.6]{Springer:intersection}, $\Delta$ is 
uniquely determined by the conditions
\begin{enumerate}
\item[(a)]
$\Delta(\tilde T_{\mathbf{s}}. m)=\tilde T_{\mathbf{s}}^{-1} . \Delta(m)$
for $\mathbf{s} \in \mathbf {S}$ and $m \in \mathcal{M}$, and
\item[(b)]
if $\tilde R_{a,b}(u)=(-u)^{-\Bruhatl(b)+\Bruhatl(a)} R_{a,b}(u^2)$
for the polynomials $R_{a,b}$ described in
\cite[\S 2]{Kazhdan:representations} or 
\cite[3.1 -- 3.7]{Dyer:Iwahori-Hecke}, then
\begin{equation*}
\tilde b_{[I,a_1,b_1],[J,a_2,1]}=\left\{
\begin{array}{ll}
0   &\mbox{ if } I \not \subseteq J \mbox{ or } b_1 \not \in W_J \\
(\bar \alpha)^{|J|-|I|} \tilde R_{a_2b_1,a_1}   &\mbox{otherwise}.
\end{array}\right.
\end{equation*}
\end{enumerate}

\section{Statement of Results} \label{section:results}
Maintain the notations and assumptions of Section
\ref{section:Springer's conjecture}. 

\subsection{}
Fix a Coxeter system
$(\hat{W},\hat{S})$ with the following properties:
\begin{enumerate}
\item[(a)] 
$W$ is the standard parabolic subgroup $W=\hat{W}_S$ of
$\hat{W}$ generated by $S \subseteq \hat{S}$ and
\item[(b)] 
there is a bijection  
$\theta\colon S\rightarrow \hat{S} \smallsetminus S$ such that for 
$r, s \in S$,
$r$ and $\theta(s)$ are joined by an edge of the Coxeter graph of
$(\hat{W},\hat{S})$ (i.e. do not commute in $\hat{W}$) iff $r=s$.
\end{enumerate}
Note that the order of $r\theta(r)$ for  $r\in S$ is assumed to be
 three or greater but is  not specified more precisely, and that  there  is no
  assumption made  about the orders of the products
 $\theta(r)\theta(s)$ for distinct $r,s\in S$. 
  In particular,  if $W$ is non-trivial, there are infinitely many
possible  choices of 
$(\hat W,\hat S)$ up to isomorphism.

For example, if $(W,S)$ is of type $A_1$ ,
then $(\hat W,\hat S)$ could be taken to be any irreducible
 dihedral Coxeter
 system, for instance   of type $A_2$.
If $(W,S)$ is of type $A_2$ (resp.,  $A_3$, $B_2$)
then $(\hat W,\hat S)$ could be chosen to be  of type  
 $A_4$ (resp.,  $E_6$, $F_4$)
with the Coxeter graph of $(W,S)$ embedded as the full subgraph
of that of $(\hat{W},\hat{S})$ on the vertex set obtained by
deleting  all 
 of the  terminal vertices of the latter graph.
 If $(W,S)$ is irreducible but not of   type
$A_1$, $A_2$, $A_3$ or $B_2$, then   $(\hat{W},\hat{S})$ cannot  be
 of finite or affine type.

\subsection{} 
Fix once and for all an enumeration of $S$ as $S=\{\,
r_1,\ldots,r_n
\,\}$ with distinct
$r_i$.  For $I \subseteq S$, define
$z_I \in \hat{W}$ as $z_I=\theta(r_{i_1})\cdots
\theta(r_{i_m})$, where $S \smallsetminus I=\{\, i_1,\ldots,
i_m \,\}$  with $i_1< \cdots <i_m$ (note $z_I$ may depend  on the
fixed enumeration of $S$).

Let $\hat{T}$ denote the set of reflections of $(\hat{W},\hat{S})$ and
set $A=\hat{T}\smallsetminus W$. 
Then by \cite[2.11]{Dyer:Hecke2}, $A$ is an initial section of  a
reflection order of $\hat{T}$. 
Correspondingly, there are a length function
$\twistl{A}\colon \hat{W} \rightarrow \mathbb{Z}$, 
a partial order $\twisto{A}$ on $\hat{W}$ and elements 
$R^A_{x,y} \in \mathcal{R}$, $P_A({x,y})\in \mathbb{Z}[q]$ defined for
$x,y \in \hat{W}$ (some of their definitions and basic properties are
summarized in Section \ref{section:twisted Chevalley-Bruhat order}.)
The following is the main result of this paper.

\begin{theorem}\label{main}
Let $\Omega=\bigcup_{I\subseteq S} Wz_IW \subseteq \hat{W}$, 
regarded as a poset in the order induced by $\twisto{A}$.
\begin{enumerate} 
\item[(i)]  
$\Omega$ is locally closed in $(\hat{W},\twisto{A})$, i.e. if
$x \twisto{A} y \twisto{A} z$ with $x,z\in \Omega$ and $y\in \hat{W}$,
then $y\in \Omega$.
\item [(ii)] 
The map $\varphi\colon V \rightarrow \Omega$ given by
$\varphi([I,a,b])=az_Ib^{-1}$ is a poset isomorphism.
\item[(iii)] 
For $v\in V$, $\twistl{A}(\varphi(v))=d(v) - |S|$.
\item[(iv)] 
$\varphi((a,b).v)=a \varphi(v) b^{-1}$ for any  $v \in V$ and 
$(a,b) \in \mathbf{W}$.
\item[(v)] 
For $v,w\in V$, $\tilde b_{w,v}=R^A_{\varphi(w),\varphi(v)}$.
\end{enumerate}
\end{theorem}

The theorem allows one to directly transfer many known properties of
the twisted Chevalley-Bruhat order $\twisto{A}$ to the poset $V$. 
For example, one obtains the following.

\begin{corollary}\label{maincor} 
Let $v,w\in V$ and $n=d(v)-d(w)$.
\begin{enumerate}
\item[(i)]
$\sum_{\substack{z \in V \\ w \leqslant z \leqslant v}} 
\tilde b_{w,z}(u) \tilde b_{z,v}(u^{-1})=\delta_{w,v}$,
where $\delta_{w,v}$ is the Kronecker delta.
\item[(ii)] 
if $w \leq v$, the closed interval $[w,v]$ and the opposite poset
$[w,v]^\op$ are pure EL-shellable posets 
{\rm(}in which every maximal chain contains $n+1$ elements{\rm)} 
in the sense of  \cite[Definition 2.1]{BjornerWachs}.
\item[(iii)] 
if $w \leq v$  and $n\geq 2$, then the order complex {\rm(}the abstract
simplicial complex with totally ordered subsets as simplexes{\rm)} of the
open interval $(w,v)$ is a combinatorial $(n-2)$--sphere.
\end{enumerate}
\end{corollary}

Corollary \ref{maincor} (i) asserts that $\Delta$ 
(or more precisely its natural extension to a map 
$\tilde{\mathcal{M}} \rightarrow \tilde{\mathcal{M}}$) is an involution;
this was proved geometrically by Springer in the case $W$ is a finite
Weyl group, and  conjectured by him to hold in general.
Once it is known, it becomes possible to define analogues $c_{w,v}$
(resp., $c_{w,v}^{\mbox{\tiny inv}}$) of the Kazhdan-Lusztig
polynomials $P_{x,y}$ (resp., inverse Kazhdan-Lusztig polynomials
$Q_{x,y}$) as in \cite{Kallen:computing}.
Regarding these, we record

\begin{corollary}\label{maincor1} 
For $w,v\in V$, we have $c_{w,v}=P_A({\varphi(v),\varphi(w)})$ and
$c_{w,v}^{\mbox{\rm\tiny inv}}=P_{\hat{T} \smallsetminus A}
({\varphi(w),\varphi(v)})$, 
where the $P_{\hat{T} \smallsetminus A}$ and $P_A$
are the polynomials to be described in 3.9.
\end{corollary}

\subsection{} 
Introduce a graph, also denoted $V$, with vertex set $V$ and an edge
from $w$ to $v$ if $\tilde b'_{w,v}(1) \ne 0$, where $p'(1)$ denotes
the value at $1$ of the derivative with respect to
$u$ of a Laurent polynomial $p\in \mathcal{R}$. 
By \cite[3.10]{Springer:intersection} (extended to general $(W,S)$)
and Theorem \ref{main}, one has an edge from $w$ to $v$ iff either
$d(w)<d(v)$ and there is a reflection $r \in \mathbf{W}$ such that
$w=rv$, or $w=[I,a,b]$, $v=[J,c,d]$, $I\subset J$, 
$| J \smallsetminus I |=1$ and there exists $f\in W_J$ with $a=cf$
and $b=df$.

\begin{corollary}\label{maincor2}
\begin{enumerate}
\item[(i)] 
If $v,w \in V$, then $w \leq v$ iff there is a directed path in the
graph $V$ from $w$ to $v$.
\item[(ii)] {\rm(Deodhar's inequality)}
For $w \leq z \leq v$ in $V$, the graph $V$ has at least
$d(v)-d(w)$ edges with $z$ as one vertex and both vertices in $[w,v]$.
\end{enumerate}
\end{corollary}

\subsection{}   
Other consequences of the theorem include results on the structure of
``quotients'' of $V$ in \cite{Dyer:quotients}  
(such as posets of shortest orbit representatives in $V$ for orbits of
parabolic subgroups of $W$ for left or right multiplication), 
the analogue for $\tilde b_{w,v}$ in terms of
the above graph for the generating function formula 
\cite[\S 3]{Dyer:Iwahori-Hecke} for $R^A_{x,y}$, 
recurrence formulas \cite[\S 4]{Dyer:Iwahori-Hecke} for the $c_{w,v}$ 
similar to the classical ones \cite{Kazhdan:representations} 
for $P_{x,y}$, existence of three-parameter versions $c_{w,v}^B$ 
of the $c_{w,v}$ with conjecturally non-negative coefficients 
(involving as an additional parameter $B$ an arbitrary initial section 
of a reflection order of $\hat W$),  
and existence of highest weight representation categories
for which the $c_{w,v}$ conjecturally describe Verma module
multiplicities etc. 
We refrain from describing any of these in detail, but record the
following.

\begin{corollary}\label{finitecor}  
Suppose that $W$ is finite with the longest element $w_S$. 
Define a map $\varphi': V \rightarrow \hat{W}$ by 
$\varphi'([I,a,b]) =a z_I b^{-1}w_S$.
\begin{enumerate} 
\item[(i)] 
The map $\varphi'$ restricts to an order isomorphism from
$V$ to the interval $[w_S z_\emptyset w_S,1]$ in reverse
Chevalley-Bruhat order.
\item[(ii)] 
$\Bruhatl(\varphi'(v))=-d(v)+| S |+\Bruhatl(w_S)$ for $v\in V$.
\item[(iii)] 
$\varphi'((a,b).v)=a \varphi'(v) w_S b^{-1} w_S$ for $v \in V$ and
$(a,b) \in \mathbf{W}$.
\item[(iv)] 
$\tilde b_{w,v}=\tilde R_{\varphi'(v),\varphi'(w)}$.
\item[(v)] 
$c_{w,v}=Q_{\varphi'(v),\varphi'(w)}$ and
$c_{w,v}^{\mbox{\rm\tiny inv}}=P_{\varphi'(v),\varphi'(w)}$.
\end{enumerate}
\end{corollary}

\begin{remark} The above corollary can be proved independently of the
theory of the orders $\leq_A$; the main step,
analogous to \ref{lemma:recursive formula for R^A} 
and proved similarly, is to show that for any
Coxeter system $(\hat{W},\hat{S})$ with finite standard parabolic 
subgroup $(W,S)$, one has
\begin{equation*} 
\tilde R_{a_1z_1w_S,a_2z_2b_2w_S}=
\begin{cases}
\tilde R_{z_1,z_2} \tilde R_{a_1b_2^{-1},a_2} 
             &\mbox{ if } b_2 \in W_{I_1} \\
0,  &\mbox{ otherwise,}
\end{cases}
\end{equation*} 
where $z_i\in \hat W_{\hat{S} \smallsetminus S}$,
$I_i=S \cap z_i S z_i^{-1}$, $a_i\in W^{I_i}$, and  $b_2 \in W$.
\end{remark}

\section{Twisted Chevalley-Bruhat Orders}
\label{section:twisted Chevalley-Bruhat order}

Let $(W,S)$ be an arbitrary Coxeter system and 
$T=\bigcup_{x \in W} xSx^{-1}$ be the set of reflections in $W$.

\subsection{Reflection orders on $T$}
By \cite[Definition 2.1]{Dyer:Hecke2}, a total order $\preceq$ on $T$ 
is called a reflection order if for any dihedral reflection subgroup
$W'$ of $W$ with its canonical simple reflections $S'=\{ r,s \}$ with 
respect to $S$, either $r \prec rsr \prec \cdots
\prec srs \prec s$ or $s \prec srs \prec \cdots \prec rsr
\prec r$.

$A \subseteq T$ is called an initial section of reflection orders on 
$T$ if there is a reflection order $\preceq$ on $T$ such that
$x \prec y$ for all $x \in A$ and $y \in T \smallsetminus A$.

\subsection{Twisted Chevalley-Bruhat order}
We give a brief description of twisted Chevalley-Bruhat orders here.
The reader can refer to \cite[\S 1]{Dyer:Iwahori-Hecke} for the details.

Regard the power set $\mathcal{P}(T)$ of $T$  as an additive abelian 
group under the symmetric difference 
$A+B=(A \cup B) \smallsetminus (A \cap B)$. 
Define a map $N\colon W \rightarrow \mathcal{P}(T)$ by
$N(x)=\{\,t \in T \mid tx \Bruhato x\,\}$ for $x \in W$.
We have $\Bruhatl(x)=|N(x)|$.
$N(x)$ may be characterized by $N(xy)=N(x) + xN(y)x^{-1}$ for  
$x,y \in W$ and $N(s)=\{\, s \,\}$ for $s \in S$.
There is an action of $W$ on $\mathcal{P}(T)$ given by 
$x.A=N(x)+xAx^{-1}$ for $x \in W$ and $A \in \mathcal{P}(T)$.

Fix an initial section $A$ of a reflection order on $T$.
Define the twisted length function $\twistl{A}$ on $W$ by
\[ \twistl{A}(x)=\Bruhatl(x)-2 |N(x^{-1}) \cap A| \qquad
\mbox{ for } x \in W. \]
The twisted Chevalley-Bruhat order $\twisto{A}$ on $W$ is generated 
by the relation $x \twisto{A} tx$ for $x \in W$ and $t \in T$ with
$\twistl{A}(x)<\twistl{A}(tx)$.
By \cite[Proposition 1.2]{Dyer:Iwahori-Hecke}, 
we have $x \twisto{A} tx$ iff $t \not\in x.A$.

\subsection{The Coxeter system $(W',S')$}
Set $S'=\{\, t \in T \mid t.A=\{t\} + A \,\}$.
Let $W'$ be the subgroup of $W$ generated by $S'$ and
set $T'=\bigcup_{x \in W'} x S' x^{-1}$.
Note that the map $N'\colon W' \rightarrow \mathcal P(T')$ given by
$N'(x)=x.A+A$ satisfies $N'(xy)=N'(x)+x N'(y) x^{-1}$ for $x,y \in W'$ 
and $N'(s)=\{\, s \,\}$ for $s \in S'$.
By \cite[\S 2]{Dyer:reflection}, $(W',S')$ is a Coxeter system.

For $x \in W$, define 
$\mathcal{L}_A(x)=\{\, s \in S \mid sx \twisto{A} x \,\}$
and $\mathcal{R}_A(x)=\{\, t \in S' \mid xt \twisto{A} x \,\}$.
By \cite[1.8]{Dyer:Iwahori-Hecke}, we see that if 
$s \in \mathcal{L}_A(x)$, then
$\twistl{A} (sx)=\twistl{A} (x)-1$ 
and that if $t \in \mathcal{R}_A(x)$, then
$\twistl{A} (xt)=\twistl{A} (x)-1$.
By \cite[Proposition 1.9]{Dyer:Iwahori-Hecke},
the twisted Chevalley-Bruhat order $\twisto{A}$ satisfies the 
following properties.

\begin{proposition} \label{proposition:Z-property}
Let $x,y \in W$, $s \in S$, and $t \in S'$.
\begin{enumerate}
\item[(i)] 
If $sx \twisto{A} x$ and $sy \twisto{A} y$, then 
$sx \twisto{A} sy$ iff $sx \twisto{A} y$ iff $x \twisto{A} y$.
\item[(ii)]
If $xt \twisto{A} x$ and $yt \twisto{A} y$, then 
$xt \twisto{A} yt$ iff $xt \twisto{A} y$ iff $x \twisto{A} y$.
\end{enumerate}
\end{proposition}

We call (i) (resp., (ii)) in Proposition \ref{proposition:Z-property}
the left (resp., right) Z-property of $\twisto{A}$.

\begin{lemma}\label{lemma:maximums} 
Regard the power set of $W$ as a monoid under the product
$XY=\{ xy\mid x\in X,y\in Y\}$.
For $w\in W$, let $[1,w]=\{ z\in W \mid 1 \Bruhato z \Bruhato w\,\}$.
Suppose that $B$ is a subset of $W$ with a maximum 
{\rm(}resp., minimum{\rm)} element in the order $\twisto{A}$.
Then $[1,w]B$ has a maximum {\rm(}resp., minimum{\rm)} element in
$\twisto{A}$.
\end{lemma}

\begin{proof} 
We treat only the case $B$ has a maximum element $m$.
It is well-known that if $r_1 \cdots r_n$ is a reduced expression 
for $w$, then $[1,w]=[1,r_1] \cdots [1,r_n]$.
Hence we may assume $w=r \in S$. Let $m'$ be the maximum of $m$ and $rm$ 
in $\twisto{A}$. 
The left $Z$-property implies that $m'$ is a maximum element of
$[1,r]B$.  
\end{proof}

\subsection{A relation on $W \times W$}
For $(x,y)$, $(x',y') \in W \times W$, 
write $(x,y) \longrightarrow' (x',y')$
(see \cite[\S 2]{Dyer:Iwahori-Hecke}) if there exists $s \in S$ 
satisfying one of the
following conditions:
\begin{enumerate}
\item[(a)]
$s \in \mathcal{L}_A (x)$, $s \in \mathcal{L}_A(y)$, 
$x'=sx$, $y'=sy$,
\item[(b)]
$s \not\in \mathcal{L}_A (x)$, $s \not\in \mathcal{L}_A(y)$, 
$x'=sx$,  $y'=sy$,
\item[(c)]
$s \not\in \mathcal{L}_A (x)$, $s \in \mathcal{L}_A(y)$, 
$x'=x$, $y'=sy$,
\item[(d)]
$s \not\in \mathcal{L}_A (x)$, $s \in \mathcal{L}_A(y)$, 
$x'=sx$, $y'=y$.
\end{enumerate}
Define the relation $\longrightarrow$ on $W \times W$ to be the 
(reflexive) transitive closure of the relation $\longrightarrow'$. 
It is known by \cite[Lemma 2.2 and Proposition 2.5]{Dyer:Iwahori-Hecke}
that if $(v,w)
\longrightarrow (1,1)$ then $v \twisto{A} w$ and the closed interval
$[v,w]$ in $\twisto{A}$ is finite.

\begin{lemma}\label{lemma:parabolic interval} 
Suppose $v,w\in W$ with $(v,w) \longrightarrow (1,1)$.
If $v \twisto{A} z\twisto{A} w$, 
$J\subseteq S$ and $wv^{-1}\in W_J$, then
$zv^{-1},wz^{-1}\in W_J$.
\end{lemma}

\begin{proof} 
This holds by \cite[Lemma 2.12]{Dyer:quotients}.
\end{proof}

\subsection{$R^A$-polynomials} \label{subsection:R-polynomial}
Recall that $\mathcal{R}=\mathbb{Z}[u,u^{-1}]$ and $\alpha=u-u^{-1}$.
By \cite[Corollary 3.6 (1)]{Dyer:Iwahori-Hecke},
there is a unique family of polynomials
$R^A_{x,y} \in \mathbb{Z}[\bar \alpha]$  defined for $x,y \in W$ with
$(x,y) \longrightarrow (1,1)$ such that
\begin{enumerate}
\item[(a)]
$R^A_{1,1}=1$;
\item[(b)]
if $x,y \in W$ satisfy $(x,y) \longrightarrow (1,1)$ and
$s \in \mathcal L_A(y)$, then
\[ R^A_{x,y}=\left\{
\begin{array}{ll}
R^A_{sx,sy}                       &\mbox{ if } s \in \mathcal L_A(x) \\
R^A_{sx,sy}+\bar\alpha R^A_{x,sy} &\mbox{ if } s \not\in \mathcal L_A(x),
\end{array} \right. \]
\end{enumerate}
where it is understood that in case $s \not\in \mathcal{L}_A(x)$ and
$sx \not\twisto{A} sy$, the term $R^A_{sx,sy}$ is zero.

\subsection{Generating function formula}
\label{subsection:generating function} 

The proof that the $R^A_{x,y}$ are
well-defined is accomplished in conjunction with the proof of a 
generating function formula
\[R^A_{x,y}=\sum_{\substack{n \\(t_1, \ldots, t_n)}}\bar{\alpha}^n,\]
for $(x,y)\longrightarrow (1,1)$, where the sum is over $n\in
\mathbb{N}$ and
$n$-tuples 
$(t_1,\ldots, t_n)$ of elements of $T$ satisfying
$x \twisto{A} t_1x \twisto{A} t_2t_1x \twisto{A} \cdots \twisto{A}
t_n \cdots t_2t_1x=y$ and 
$t_1 \preceq  t_2\preceq \cdots \preceq t_n$ with $\preceq$ 
a fixed but arbitrary reflection order on $T$.

By \cite[Corollaries 3.6 and 3.7]{Dyer:Iwahori-Hecke},  
if $(x,y) \longrightarrow (1,1)$,
then 
$\sum_{x \twisto{A} z \twisto{A} y} R^A_{x,z}(u)R^A_{z,y}(u^{-1})
=\delta_{x,y}$ and
$R^A_{x,y}$ is a monic polynomial in $\bar\alpha$ of degree
$\twistl{A}(y)-\twistl{A}(x)$. 
We also note that if $A=\emptyset$, then
$R^A_{x,y}=\tilde R_{x,y}$ is exactly the polynomial described in 
\ref{subsection:the map Delta}.

\subsection{}  From \cite[\S 4]{Dyer:Iwahori-Hecke},
there are
unique
$p_A(v,w)\in
\mathbb{Z}[u,u^{-1}]$ defined for $(v,w) \longrightarrow (1,1)$, 
such that
\begin{enumerate}\item[(a)] $p_A(w,w)=1$
\item[(b)] $p_A(v,w)\in u^{-1}\mathbb{Z}[u^{-1}]$ if $v\neq
w$
\item[(c)]$p_A(v,w)
=\sum_{v\leq_A z \leq_A w} R_A(v,z)\overline{p_A(z,w)}$ 
\end{enumerate}
We define the ``Kazhdan-Lusztig'' polynomials   
$P_A(v,w)=P_A(v,w)(q)\in\mathbb{Z}[q]$ for $\leq_A$ 
 by the formula
$u^{l_A(w)-l_A(v)}p_A(v,w)=P_A(v,w)(u^2)$, if
$(v,w)\longrightarrow (1,1)$.

The following lemma summarizes some simple relations
we shall require between objects associated to varying initial
sections $A$.

\begin{lemma}\label{varyingA} 
Let $A$ be an initial section of reflection orders on $T$ and $v,w\in W$.
Then
\begin{enumerate}
\item[(i)]  
$T+A$ and $x.A$ for $x\in W$ are initial sections of reflection
orders on $T$.
\item[(ii)] 
We have $\twistl{T+A}(w)=-\twistl{A}(w)$ and 
$(W,\twisto{T+A})=(W,\twisto{A})^{\op}$, 
$R^{T+A}_{v,w}=R^A_{w,v}$.
\item[(iii)] 
For a fixed $x\in W$, the map $w \mapsto  wx$ defines a
poset isomorphism 
$\phi\colon (W,\twisto{x. A}) \rightarrow (W,\twisto{A})$ satisfying
$\twistl{A}(\phi(w))=\twistl{x.A}(w)+\twistl{A}(x)$,
$R^{A}_{\phi(v),\phi(w)}=R^{x.A}_{v,w}$, and
$P_{A}(\phi(v),\phi(w))=P_{x.A}({v,w})$.
\end{enumerate}
{\rm(}here, equality of $R^A$-polynomials or $P_A$-polynomials is 
interpreted in the sense
that if one side is not defined then neither is the other.{\rm)}
\end{lemma}

\begin{proof} Part (i) is in \cite[Lemma 2.7]{Dyer:Hecke2}. 
For (ii) and (iii),
it is enough by the definitions of $\twisto{A}$, $R^A$, $P_A$ to prove 
the assertions on the length functions; for (ii), this is immediate
from the definition, and for (iii), it is established in the proof of
\cite[Proposition 1.1]{Dyer:Iwahori-Hecke}.
\end{proof}

\subsection{Modules over the Iwahori-Hecke algebra}

Let $\mathcal{H}$ and $\mathcal{H}'$ be the Iwahori-Hecke algebras of
$(W,S)$ and $(W',S')$, respectively.
Denote the standard basis of $\mathcal{H}'$ over $\mathcal{R}$ by
$\{\, \tilde T_x' \mid x \in W \,\}$.
Let $\mathcal{H}_A$ be the set of formal $\mathcal{R}$-linear 
combinations $\sum_{x \in W} a_x \tilde t_x$ such that there exists 
$y \in W$ so that $a_x=0$ unless $x \twisto{A} y$.
$\mathcal{H}_A$ admits an $(\mathcal{H},\mathcal{H}')$-bimodule 
structure (see \cite[\S 4]{Dyer:Iwahori-Hecke}).
If $s \in S$, 
then $\tilde T_s (\sum_{x \in W} a_x \tilde t_x)=\sum_{x \in W} b_x
\tilde t_x$, where
\[ b_x=\left\{
\begin{array}{ll}
a_{sx}               &\mbox{ if } s \not\in \mathcal{L}_A(x) \\
a_{sx}+\alpha a_x    &\mbox{ if } s \in \mathcal{L}_A(x).
\end{array} \right. \]
If $t \in S'$, then 
$(\sum_{x \in W} a_x \tilde t_x)\tilde T_t'=\sum_{x \in W} c_x 
\tilde t_x$, where
\[ c_x=\left\{
\begin{array}{ll}
a_{xt}               &\mbox{ if } t \not\in \mathcal{R}_A(x) \\
a_{xt}+\alpha a_x    &\mbox{ if } t \in \mathcal{R}_A(x).
\end{array} \right. \]

In the rest of this section, we  assume for simplicity that 
$x \twisto{A} y$ in $W$  iff $(x,y) \longrightarrow (1,1)$.

\subsection{An Involution on $\mathcal{H}_A$}
By \cite[\S 1]{Kazhdan:representations}, there are semi-linear ring 
involutions on $\mathcal{H}$ and  $\mathcal{H}'$ given by
\[\overline{\sum_{x \in W} a_x \tilde T_x}
=\sum_{x \in W} \bar a_x (\tilde T_{x^{-1}})^{-1},\qquad
\overline{\sum_{x \in W} a_x \tilde  T_x'}
=\sum_{x \in W} \bar a_x (\tilde T_{x^{-1}}')^{-1}\qquad . \]
(Here, semi-linearity is with respect to the involution 
$a \mapsto \bar{a}$ of $\mathcal{R}$.)
By \cite[Proposition 4.2 (1)]{Dyer:Iwahori-Hecke}, 
the $\mathcal R$-module $\mathcal{H}_A$ admits a semi-linear involution
defined by
\[ \overline{\sum_{x \in W} a_x \tilde t_x}
=\sum_{\substack{x,y \in W \\ x \twisto{A} y}} \bar a_{y} R^A_{x,y} (u)
\tilde t_x. \]

\begin{lemma}\label{lemma:compatible} 
These involutions are compatible; namely, 
$\overline{h \cdot m \cdot h'}=\overline{h}\cdot\overline{m}\cdot
\overline{h'}$ for
$m \in \mathcal{H}_A$, $h\in \mathcal{H}$, and $h' \in \mathcal{H}'$.
\end{lemma}

\begin{proof} 
See \cite[Proposition 4.2 (2) and Corollary 4.16]{Dyer:Iwahori-Hecke}.
\end{proof}

The following is the right analogue of
\ref{subsection:R-polynomial} (b).

\begin{lemma} \label{lemma:recursive formula}
If $x,y \in W$ and $t \in S'$ satisfy $(x,y) \longrightarrow (1,1)$ and
$t \in \mathcal R_A(y)$, then
\[ R^A_{x,y}=\left\{
\begin{array}{ll}
R^A_{xt,yt}                &\mbox{ if } t \in \mathcal{R}_A(x) \\
R^A_{xt,yt}+\bar\alpha R^A_{x,yt} 
                           &\mbox{ if } t \not\in \mathcal{R}_A(x).
\end{array} \right. \]
\end{lemma}

\begin{proof}
It is easy to compute
\[ \overline{(\tilde t_x \tilde T_t')}=\left\{
\begin{array}{ll}
\sum_{z \in W} (R^A_{xt,z}+\bar\alpha R^A_{x,z}) \tilde t_z 
         &\mbox{ if } t \in \mathcal{R}_A(x) \\
\sum_{z \in W} R^A_{xt,z} \tilde t_z   
         &\mbox{ if } t \not\in \mathcal{R}_A(x)
\end{array} \right. \]
\begin{align*}
   \overline{(\tilde t_x)} \, \overline{(\tilde T_t')}
&=\left( \sum_{z \in W} R^A_{x,z} \tilde t_z \right)
   \left( \tilde T_t'-\alpha \tilde T_1' \right)          \\
&=\sum_{\substack{z \in W \\ t \in \mathcal R_A(z)}} R^A_{x,zt} 
    \tilde t_z
   +\sum_{\substack{z \in W \\ t \not\in \mathcal R_A(z)}}
    (R^A_{x,zt}+\bar\alpha R^A_{x,z} ) \tilde t_z.
\end{align*}
Since $\overline{(\tilde t_x  \tilde T_t')}
=\overline{(\tilde t_x)} \, \overline{(\tilde T_t')}$,
comparing the coefficients of $\tilde t_{yt}$ gives the desired 
recurrence formula.
\end{proof}

\section{Twisted Chevalley-Bruhat Order Associated to Standard 
Parabolic Subgroups}
\label{section:special order}
In this section, $(\hat{W},\hat{S})$ denotes an arbitrary Coxeter
system. 
Fix $S \subseteq \hat{S}$ and set $W=\hat{W}_S$.
By \cite[Proposition 2.11]{Dyer:Iwahori-Hecke}, 
$A:=\hat{T}\smallsetminus W$ is an initial section of a reflection order
of $\hat{W}$ and the corresponding partial order $\twisto{A}$ on $\hat{W}$
satisfies
$v \twisto{A} w$ iff $(v,w) \longrightarrow (1,1)$.

\begin{lemma}\label{lemma:projection} 
Let either $B=A$ or $B=\emptyset$.
\begin{enumerate}
\item[(i)] 
$S\subseteq \{\, t\in \hat{T}\mid t \cdot B=B+\{t\} \,\}$.
\item[(ii)] For $w\in\hat{W}^S$ and $x\in W$,
$\twistl{A}(wx)=\Bruhatl(x)-\Bruhatl(w)$ and
$\Bruhatl (wx)=\Bruhatl (x)+\Bruhatl (w)$.
\item[(iii)] 
Every element $z$ of $\hat{W}^S$ is the minimum element of its coset
$zW$ in the order $\twisto{B}$.
\item[(iv)] 
The map $\pi\colon \hat{W} \rightarrow \hat{W}^S$ such that
$\pi(w)\in \hat{W}^S\cap w W$ is order-preserving for $\twisto{B}$,
i.e. $x \twisto{B} y$ implies $\pi(x)\twisto{B} \pi(y)$.
\item[(v)] 
The restrictions of $\twisto{A}$ and reverse Chevalley-Bruhat order to
partial orders on $\hat{W}^S$ coincide.
\end{enumerate}
\end{lemma}

\begin{proof} 
Part (i) follows by simple calculations which we omit.
For (ii) and (v), see \cite[Lemma 5.3 (1) and (3)]{Dyer:Iwahori-Hecke}.
For (iii), note that if $x\in W$ has a reduced expression
$s_1 \cdots s_m$, then for $z \in \hat{W}^S$, 
$z \twisto{B} zs_1 \twisto{B} zs_1s_2 \twisto{B} \cdots \twisto{B} zx$
by (ii) and the definition of $\twisto{B}$. 
Finally, for (iv),  suppose in the proof of (iii) that 
$w \twisto{B} zx$.
Repeated application of the right $Z$-property for
$\twisto{B}$ gives some $v \in W$ with $wv \twisto{B} z$.
Then by (iii), $\pi(w)=\pi(wv)\twisto{B} wv\twisto{B} z=\pi(zx)$.
\end{proof}

\begin{lemma} \label{lemma:double coset expression}
Let $z \in \hat W_{\hat S \smallsetminus S}$ and $I_z=S\cap zSz^{-1}$.

\begin{enumerate}
\item[(i)]
$I_z=\{\, s\in S\mid sr=rs \mbox{ for all } r\in \supp (z) \,\}$,
where for any $x \in \hat{W}$, $\supp (x)$ stands for the set of all
the simple reflections occurring in some {\rm(}equivalently, all{\rm)} 
reduced expression of $x$.
\item[(ii)]
$z$ is the unique element of minimal length in the double coset $WzW$.
\item[(iii)]
For $a\in W^{I_z}$ and $b\in W$ we have
$\Bruhatl(a z b)=\Bruhatl(a) + \Bruhatl (z) + \Bruhatl (b)$; 
in particular, $az \in \hat{W}^S$.
\item[(iv)]
Every element of $WzW$ is uniquely expressible in the form
$a z b$, where $a \in W^{I_z}$ and $b \in W$.
\item[(v)] 
$W\cap zWz^{-1}=W_{I_z}$.
\item[(vi)] 
For $a\in W^{I_z}$ and $b\in W$,
$\twistl{A} (azb)=-\Bruhatl (a) -\Bruhatl (z) +\Bruhatl (b)$.
\item[(vii)] ${\hat W}_{\hat S\setminus S}\subseteq \hat{W}^S\cap
(\hat{W}^S)^{-1}$. 
\end{enumerate}
\end{lemma}

\begin{proof}
Clearly, $\Bruhatl (zr)=\Bruhatl (rz)=\Bruhatl (z)+1$ for all 
$r\in S$, which implies (vii).
  It is well-known that this implies
(ii)--(v) (see \cite[Proposition 2.7.5]{Carter:finite}). 
Then (vi) follows from (iii) and \ref{lemma:projection} (ii).

(i) Clearly the right hand side is contained in the left
hand side. For the reverse inclusion, it will suffice to
show that if
$z, z'
\in
\hat{W}_{\hat{S}\backslash S}$, $s,t\in S$ and $sz = z't$ then
$z = z'$, $s=t$ and $s$ commutes with all $r \in \supp(z)$.
To see this, note that since $\supp(sz)=\supp(z) \cup \{s\}$ 
and $\supp(zt)=\supp(z) \cup \{t\}$ are disjoint unions, 
we must have $s=t$. 
 The remainder of the proof is by induction on $\Bruhatl (z)$. Write
$z = s_{1}\ldots s_{m}$ (reduced) and let $m\geq 1$. By the
exchange property
$ss_{1}\ldots s_{m-1}$ equals either $z'$ or $z''s$ with
$z''\leq _{\emptyset}z'$. The first case is impossible,
because $s\not \in \rm{Supp}(z')$. By induction $z'' =
s_{1}\ldots s_{ m-1}$, and $s$ commutes with $s_{1}, \ldots,
s_{m-1}$. Then \[ss_{1}\ldots s_{m} = s_{1}\ldots
s_{m-1}ss_{m} = s_{1}\ldots s_{m}s,\] and one sees that $s$
also commutes with $s_{m}$.  
\end{proof}

\begin{lemma} \label{lemma:middle elements in the twisted 
Chevalley-Bruhat order}
Assume that $z_1,z_2 \in \hat{W}_{\hat S \smallsetminus S}$,
$a_1 \in W^{I_{z_1}}$, $a_2 \in W^{I_{z_2}}$, and $b_1,b_2 \in W$ 
satisfy $a_1z_1b_1 \twisto{A} a_2z_2b_2$. 

\begin{enumerate}
\item[(i)]
If $x \in \hat W$ satisfies 
$a_1z_1b_1 \twisto{A} x \twisto{A} a_2z_2b_2$, 
then there exist $z_3 \in \hat W_{\hat S \smallsetminus S}$,
$a_3 \in W^{I_{z_3}}$, and $b_3 \in W$ such that $x=a_3z_3b_3$.
\item[(ii)] 
$z_2 \Bruhato z_1$.
\end{enumerate}
\end{lemma}

\begin{proof}
We use the map $\pi$ from Lemma \ref{lemma:projection}.
Note $\pi(a_iz_ib_i)=a_iz_i$ and 
$\pi (z_{i}^{-1}a_{i}^{-1})=z_{i}^{-1}$ by 
Lemma \ref{lemma:double coset expression}. 
We get by Lemma \ref{lemma:projection} (iv) that
$a_1z_1\twisto{A} \pi(x)\twisto{A} a_2z_2$; 
so by Lemma \ref{lemma:projection} (v),
$a_2z_2 \Bruhato  \pi(x) \Bruhato a_1z_1$. 
Taking inverses throughout, applying $\pi$ then taking inverses again
gives
$z_2\Bruhato v:=
\Bigl(\pi\bigl(\pi(x)^{-1}\bigr)\Bigr)^{-1}\Bruhato z_1$.  
Since $v$ is a subword of $z_1$ (or by Lemma
\ref{lemma:parabolic interval}) we have
$v\in \hat {W}_{\hat{S}\smallsetminus S}$; 
by definition, $v \in W x W$. 
Applying Lemma \ref{lemma:double coset expression} gives (i), 
and we get (ii) by taking $x=a_1z_1b_1$ above.
\end{proof}

\begin{lemma}\label{lemma:first recursive formula for R^A}
Let $z_1,z_2 \in \hat{W}_{\hat{S} \smallsetminus S}$. 
Then $R^A_{z_1,z_2}=\tilde R_{z_2,z_1}$.
\end{lemma}

\begin{proof}
Observe that for any $z \in \hat W_{\hat S \smallsetminus S}$ and
$s \in \hat S \smallsetminus S$, we have $z \twisto{A} sz$ iff
$sz \Bruhato z$.
We shall proceed by induction on $\Bruhatl(z_1)$.

If $\Bruhatl(z_1)=0$, then $z_1=1$. 
If either $R^A_{z_1,z_2}$ or $\tilde R_{z_2,z_1}$ is
non-zero, then by Lemma \ref{lemma:projection}(v) and
\ref{subsection:generating  function} it follows that $z_2=1$
and that
$R^A_{z_1,z_2}=1=\tilde R_{z_2,z_1}$.

Suppose that $R^A_{z_1',z_2'}=\tilde R_{z_2',z_1'}$
for any $z_1',z_2' \in \hat W_{\hat{S} \smallsetminus S}$ with 
$\Bruhatl(z_1') < \Bruhatl(z_1)$.
Pick an $s \in \hat{S}$ such that $sz_1 \Bruhato z_1$.
If $sz_2 \Bruhato z_2$, then by the induction hypothesis and
the recurrence formula in \ref{subsection:R-polynomial}, we obtain
$R^A_{z_1,z_2}=R^A_{sz_1,sz_2}=\tilde R_{sz_2,sz_1}=\tilde R_{z_2,z_1}$.
If $z_2 \Bruhato sz_2$, then
$ R^A_{z_1,z_2}
=R^A_{sz_1,sz_2}+\bar \alpha R^A_{z_1,sz_2}
=R^A_{sz_1,sz_2}+\bar \alpha R^A_{sz_1,z_2}
=\tilde R_{sz_2,sz_1}+\bar \alpha \tilde R_{z_2,sz_1}
=\tilde R_{z_2,z_1}. $
\end{proof}

\begin{lemma} \label{lemma:recursive formula for R^A}
Assume that  $z_1,z_2 \in \hat W_{\hat{S} \smallsetminus S}$,
$a_1 \in W^{I_{z_1}}$, $a_2 \in W^{I_{z_2}}$, and $b_1 \in W$.
Then
\[ R^A_{a_1z_1b_1,a_2z_2}=
\tilde R_{z_2,z_1} \tilde R_{a_2 b_1^{-1},a_1}\chi(I_{z_2},b_1),\]
where 
$\chi(I,b):=\begin{cases}
1 &\mbox{ if } b\in \hat{W}_I \\ 
0 &\mbox{ otherwise.}
\end{cases}$
\end{lemma}

\begin{proof}
We shall proceed by induction on $\Bruhatl (a_1)$.

If $\Bruhatl(a_1)=0$, then $a_1=1$.
If the LHS is non-zero, 
then $z_1\twisto{A} z_1b_1 \twisto{A} a_2z_2\twisto{A} z_2$, 
so $b_1=a_2=1$ by Lemma \ref{lemma:parabolic interval}. 
If the RHS is non-zero,
then $a_2b_1^{-1}\Bruhato 1$ implies $a_2=b_1=1$ since
$\Bruhatl(a_2b_1^{-1})=\Bruhatl(a_2)+\Bruhatl(b_1)$. 
Hence if either side is non-zero, then
$a_2=b_1=1$ and the assertion reduces to 
Lemma \ref{lemma:first recursive formula for R^A}.
On the other hand, if both sides are zero, they are equal.

Suppose that
$R^A_{a_1'z_1'b_1',a_2'z_2'}=\tilde R_{z_2',z_1'}
\tilde R_{a_2'(b_1')^{-1},a_1'}\chi(I_{z'_2},b'_1)$
whenever $z_1',z_2' \in \hat W_{\hat{S} \smallsetminus S}$,
$a_1' \in W^{I_{z_1'}}$, $a_2' \in W^{I_{z_2'}}$, and $b_1' \in W$
subject to $\Bruhatl (a_1') < \Bruhatl (a_1)$.
Pick an $s \in S$ with $sa_1 \Bruhato a_1$.
By Lemma \ref{lemma:double coset expression} (vi), we have
$a_1z_1b_1 \twisto{A} sa_1z_1b_1$
(we omit further reference to our use of  Lemma
\ref{lemma:double coset expression} (vi) below in view of
their  frequency). 
We need to deal with three different situations.

(1) If $s a_2 \Bruhato a_2$, then $a_2 z_2 \twisto{A} s a_2 z_2$.
We obtain by \ref{subsection:R-polynomial} and induction that
\begin{align*}
R^A_{a_1 z_1 b_1, a_2 z_2}
&=R^A_{sa_1 z_1 b_1, s a_2 z_2}  \\
&=\tilde R_{z_2,z_1} \tilde R_{sa_2 b_1^{-1}, sa_1}\chi(I_{z_2},b_1)  \\
&=\tilde R_{z_2,z_1} \tilde R_{a_2 b_1^{-1}, a_1}\chi(I_{z_2},b_1).
\end{align*}

(2) If $a_2 \Bruhato sa_2$ and $sa_2 \in W^{I_{z_2}}$, then
$sa_2 z_2 \twisto{A} a_2z_2$.
We obtain
\begin{align*}
R^A_{a_1z_1b_1, a_2z_2}
&=R^A_{sa_1z_1b_1, sa_2z_2}+\bar\alpha R^A_{a_1z_1b_1,sa_2z_2} \\
&=R^A_{sa_1z_1b_1, sa_2z_2}+\bar\alpha R^A_{sa_1z_1b_1,a_2z_2} \\
&=(\tilde R_{z_2,z_1} \tilde R_{sa_2b_1^{-1},sa_1}
   +\bar\alpha \tilde R_{z_2,z_1}\tilde R_{a_2 
b_1^{-1},sa_1})\chi(I_{z_2},b_1) \\
&=\tilde R_{z_2,z_1} \tilde R_{a_2b_1^{-1},a_1}\chi(I_{z_2},b_1).
\end{align*}

(3) If $a_2 \Bruhato sa_2$ and $sa_2=a_2 t$ for some $t \in W_{I_{z_2}}$,
then $a_2 z_2 \twisto{A} a_2 z_2 t=s a_2 z_2$. 
Note $\chi(I_{z_2},b_1t)=\chi(I_{z_2},b_1)$.
We obtain
$R^A_{a_1z_1b_1,a_2z_2}
=R^A_{sa_1 z_1 b_1, sa_2z_2}
=R^A_{sa_1 z_1 b_1, a_2z_2t}.$
If $b_1t \Bruhato b_1$, then Lemma \ref{lemma:recursive formula} gives
\begin{align*}
R^A_{a_1z_1b_1,a_2z_2}
&=R^{A}_{sa_1 z_1 b_1t,a_2z_2} \\
&=\tilde R_{z_2,z_1} \tilde R_{a_2 t b_1^{-1},s a_1}\chi(I_{z_2},b_1) \\
&=\tilde R_{z_2,z_1} \tilde R_{s a_2 b_1^{-1}, s a_1}\chi(I_{z_2},b_1) \\
&=\tilde R_{z_2,z_1} \tilde R_{a_2 b_1^{-1},a_1}\chi(I_{z_2},b_1).
\end{align*}
If $b_1 \Bruhato b_1t$, then using Lemma \ref{lemma:recursive formula}
again,
\begin{align*}
R^A_{a_1z_1b_1,a_2z_2}
&=R^A_{sa_1z_1b_1,a_2z_2}+\bar\alpha R^A_{s a_1z_1b_1 t, a_2z_2} \\
&=(\tilde R_{z_2,z_1} \tilde R_{a_2 b_1^{-1},sa_1}+\bar\alpha
   \tilde R_{z_2,z_1} \tilde R_{a_2 t b_1^{-1},sa_1})\chi(I_{z_2},b_1) \\
&=\tilde R_{z_2,z_1} (\tilde R_{a_2 b_1^{-1},s a_1}+\bar\alpha
   \tilde R_{s a_2 b_1^{-1},s a_1}  )\chi(I_{z_2},b_1)\\
&=\tilde R_{z_2,z_1} \tilde R_{a_2 b_1^{-1}, a_1}\chi(I_{z_2},b_1).
\end{align*}
\end{proof}

\section{Proof of Results} \label{section:proofs}
In this section,
we fix a Coxeter system $(W,S)$ with finitely many simple reflections
and let $(\hat{W},\hat{S})$, 
$z_I\in \hat{W}_{\hat{S} \smallsetminus S}$ for 
$I\subseteq S$, $\Omega$, $\varphi$ be as in Section
\ref{section:results}. 
We let $A$, $I_z$ for $z\in \hat{W}_{\hat{S}\smallsetminus S}$,
$\twisto{A}$, $R^A_{x,y}$ etc be associated to $(\hat{W},\hat{S})$ and
$(W,S)$ as in Section \ref{section:twisted Chevalley-Bruhat order}.

\subsection{Proof of Theorem \ref{main}} 
Note that $z_\emptyset$ is a Coxeter element of
$\hat W_{\hat{S}\smallsetminus S}$ and for $I,J \subseteq S$,
$z_I\twisto{A} z_J$ iff $z_J \Bruhato z_I$ iff 
$I\subseteq J$, using \ref{lemma:projection}(v) and the
definition of the $z_I$.  From Lemma \ref{lemma:double coset
expression} (i), one sees that
$I_{z_J}=J$ for $J\subseteq S$. 
Now Theorem \ref{main} (i) holds by Lemma
\ref{lemma:middle elements in the twisted Chevalley-Bruhat order}. 
By Lemma \ref{lemma:double coset expression} (iv) and (vi), $\varphi$ is
a bijection and Theorem \ref{main} (iii) holds.
The definitions  give Theorem \ref{main} (iv) first for $(a,b)=(r,1)$ or
$(a,b)=(1,r)$ with $r\in S$ and then it follows immediately in general.
Let
$\mathcal{H}'_A$ denote the $\mathcal{R}$-submodule of $\mathcal{H}_A$
with a basis $\tilde t_x$ for $x \in \Omega$. 
It becomes an $\mathbf{H}$-module under the $\mathcal{R}$-linear action
given by
$\tilde T_{(a,b)} m=\tilde T_a m \tilde T_{b^{-1}}'$ and the map
$\beta\colon \mathcal{M}\rightarrow \mathcal{H}'_A$ given by 
$\tilde m_v \mapsto \tilde t_{\varphi(v)}$ is immediately seen to be an
$\mathbf{H}$-module isomorphism from Theorem \ref{main} (iii) and (iv), 
since $\varphi$ is bijective. 
This isomorphism extends to an $\mathbf{H}$-module isomorphism
$\hat{\beta}: \hat{\mathcal{M}} \rightarrow
\hat{\mathcal{H}}_A'$ of their completions 
(consisting of formal $\mathcal{R}$-linear
combinations of their standard basis elements).

Define $\Delta'\colon \mathcal{M}\rightarrow \hat{\mathcal{M}}$ 
as $\Delta'(m)=\hat{\beta}^{-1}(\overline{\beta  m})$,
where
$\overline{\sum_{x \in \Omega}{a_x\tilde t_x}}
=\sum_{x,y \in \Omega}
\overline{a_x}R^A_{y,x}\tilde t_y$. 
That is,
$\Delta'(\tilde m_v)=\sum_{w\in V} R^A_{\varphi(w),\varphi(v)}
\tilde m_w$  for $v \in V$. 
It is easy to see from \ref{subsection:R-polynomial} with
$A=\emptyset$ that for $I,J\subseteq S$, $\tilde R_{z_I,z_J}$
is equal to $(\bar{\alpha})^{|J|-|I|}$ if $I\subseteq J$ and is zero
otherwise. By Lemma
\ref{lemma:recursive formula for R^A} and an argument similar
to the proof of Lemma
\ref{lemma:recursive formula},  it easily follows that
$\Delta'$ has the properties
\ref{subsection:the map Delta} (a) and (b), so
$\Delta'=\Delta$, establishing Theorem \ref{main} (v).

Finally, the fact that $\varphi$ is a poset isomorphism follows from
(v), since $v \leq w$ in $V$ iff $b_{w,v}\neq 0$ and
$x\twisto{A} y$ in $\Omega$ iff $R^A_{x,y}\neq 0$ from
\ref{subsection:generating function}.

\subsection{Proofs of the corollaries}
Corollary \ref{maincor} (i) follows immediately from Theorem \ref{main}
and
\ref{subsection:generating function}.
The proofs of Corollary \ref{maincor} (ii) and (iii) are similar
using  \cite[Proposition 3.9 and Corollary 3.10]{Dyer:Iwahori-Hecke}.

Corollary \ref{maincor1}
follows immediately from the definitions and Theorem \ref{main}.

For Corollary \ref{maincor2} (i), 
note $b_{w,v}'(1)\neq 0$ is equivalent to
$(R^A_{\varphi(w),\varphi(v)})'(1)\neq 0$ 
which is equivalent in turn to the requirement that as a polynomial in
$\bar{\alpha}$, the coefficient of
$\bar{\alpha}$ in $R^A_{\varphi(w),\varphi(v)}$ is non-zero.
From \ref{subsection:generating function}, this is in turn equivalent to
the condition that
$\varphi(w)=t\varphi(v) \twisto{A} \varphi(v)$ for some $t\in \hat{T}$.
Now Corollary \ref{maincor2} (i) and (ii) follow from the definition 
of $\twisto{A}$ and \cite{Dyer:nilHecke} respectively.

For the proof of Corollary \ref{finitecor}, assume $W$ is finite. 
We have 
\[w_S.A=N(w_S)+w_SAw_S=T+w_S(\hat{T}\smallsetminus W)w_S
=T+(\hat{T} \smallsetminus T)=\hat{T},\] 
so $\twisto{w_S. A}$ is the reverse Chevalley-Bruhat order. 
It now follows from Theorem \ref{main} and Lemma \ref{varyingA} that
$\varphi'$ defines an order isomorphism of $V$ with a locally closed 
subset $\Omega$ of $W$ in reverse Chevalley-Bruhat order, satisfying all
listed properties except perhaps $\Omega=[w_Sz_\emptyset w_S,1]$. 
But $\Omega=[1,w_S][1,z_\emptyset][1,w_S]w_S$ has
maximum and minimum elements in reverse Chevalley-Bruhat order by Lemma
\ref{lemma:maximums}, 
so $\Omega$ is an interval since it is locally closed. 
We have $\Bruhatl (1)=0$ 
(resp., $\Bruhatl (w_Sz_\emptyset w_S)=2\Bruhatl(w_S)+|S|$) 
which is obviously minimal (resp., maximal) among the values 
$\Bruhatl(c)$ for $c \in \Omega$, 
so $\Omega=[w_Sz_\emptyset w_S,1]$ as required.


\end{document}